\documentclass{article}
\usepackage{amsmath,amssymb,amsthm,latexsym}

\newtheorem{thm}{Theorem}

\newtheorem{lem}[thm]{Lemma}
\newtheorem{prop}[thm]{Proposition}

\newcommand{\C}{\mathbb C}
\newcommand{\R}{\mathbb R}

\DeclareMathOperator{\Prob}{Prob}
\DeclareMathOperator{\Exp}{Exp}
\DeclareMathOperator{\Var}{Var}
\DeclareMathOperator{\tr}{tr}
\DeclareMathOperator{\K}{K}
\DeclareMathOperator{\T}{T}

\begin{document}

\title
{\bf On the distribution of the length of the second row of
a Young diagram under Plancherel measure }
\author{{\bf Jinho Baik}\\Courant
Institute of Mathematical Sciences, New York\\baik@cims.nyu.edu\\ 
{\bf Percy Deift}\\Courant Institute of Mathematical Sciences, New York\\
and \\ Institute for Advanced Study, New Jersey 
\\deift@cims.nyu.edu\\
{\bf Kurt Johansson}\\ Royal Institute of Technology, Sweden\\
kurtj@math.kth.se}
\date{January 15, 1999}
\maketitle

\begin{abstract}
 We investigate the probability distribution of the length of 
the second row of a Young diagram of size $N$ equipped with 
Plancherel measure.
We obtain an expression for the  generating function of 
the distribution in terms of 
a derivative of an associated Fredholm determinant, 
which can then be used to show that as $N\to\infty$ the distribution 
converges to the Tracy-Widom distribution [TW] for the second 
largest eigenvalue of a random GUE matrix.
This paper is a sequel to [BDJ], where we showed that as $N\to\infty$ 
the distribution of the length of the first row of a Young diagram, 
or equivalently, the length of the longest increasing subsequence 
of a random permutation, converges to the Tracy-Widom distribution [TW]
for the largest eigenvalue of a random GUE matrix.
\end{abstract}


\section{Introduction}\label{sec:intro}
Let $Y_N$ denote the set of all Young diagrams of size $N$.
As is well known, 
$Y_N$ may be viewed, equivalently, 
as the set of all partitions of $N$.
We write $\mu=(\mu_1,\mu_2,\cdots)\in Y_N$ and/or  
$\mu\vdash N$, depending on the context. 
General references for properties of the Young diagrams are, 
for example, \cite{Sa} and Section 5.1.4 in \cite{Kn}. 
For $\mu\in Y_N$, $d_\mu$ is defined to be the number of 
(standard) Young tableaux of shape $\mu$.
By Plancherel measure on $Y_N$, we mean the probability measure
defined by
\begin{equation}\label{eq1}
   \Prob\bigl(\mu\bigr) := \frac{d_\mu^2}{N!}.
\end{equation}
Under the Robinson-Schensted correspondence between $S_N$, 
the group of permutation of $\{1,2,\cdots, N\}$, and the set of pairs
of the Young tableaux of the same shape, 
Plancherel measure \eqref{eq1} is just the push forward 
of the uniform probability distribution on $S_N$.
Let $l_N^{(k)}(\mu)$ denote the number of boxes in the $k^{\text{th}}$ row
of $\mu\in Y_N$.
Schensted \cite{Sc} showed that the length of the 
longest increasing subsequence of $\pi\in S_N$ is $l_N^{(1)}(\mu)$ 
where $\mu$ is the diagram of either of the Young tableaux corresponding  
to $\pi$ under the Robinson-Schensted correspondence.
Furthermore, Greene \cite{Gr} obtained a combinatorial 
interpretations of $l_N^{(k)}$ for general $k\geq 2$, 
namely $l_N^{(1)}(\mu)+l_N^{(2)}(\mu)+ \cdots +l_N^{(k)}(\mu)$ 
is the length of the longest 
$k$-increasing subsequence of $\pi$. 
(Recall that a $k$-increasing subsequence is a union 
of $k$ disjoint increasing 
subsequences in $\pi$.)
Similarly the sum of the lengths of the first $k$ columns 
gives the length of the longest $k$-decreasing subsequence.


For $N\geq 1, n\geq 0$, define
\begin{equation}\label{e-1}
   q^{(k)}_{n,N}:= \Prob\bigl(l^{(k)}_N\le n\bigr)
=\frac1{N!}\sum_{\substack{\mu\vdash N \\ l^{(k)}_N(\mu)\le n }} d_\mu^2.
\end{equation}
Set $q^{(k)}_{n,0}:=1$, for $n\geq 0$.
Define the exponential generating function (or Poissonization) 
of $q^{(k)}_{n,N}$ by 
\begin{equation}\label{e-2}
   \phi^{(k)}_n(\lambda) := \sum_{N=0}^{\infty}
\frac{e^{-\lambda}\lambda^N}{N!} q^{(k)}_{n,N}. 
\end{equation}

First, we summarize the results of \cite{BDJ} 
using the above notation.
For $0<t\le 1$, let $u(x;t)$ be the unique solution of the
Painlev\'e II equation
\begin{equation}\label{eq4}
  u_{xx}=2u^3+xu, 
\end{equation}
with the boundary condition 
\begin{equation}\label{eq4'}
u(x;t)\sim -\sqrt{t}Ai(x)\quad\text{as}\quad x\to +\infty.
\end{equation}
The proof of the existence and the uniqueness of this solution, 
as well as the asymptotics as $x\to -\infty$, can be found, 
for example, in \cite{DZ2}.
Define the Tracy-Widom \cite{TW} distributions 
\begin{equation}\label{e-Ft}
  F(x;t) = \exp \biggl(-\int_x^\infty (y-x)(u(y;t))^2 dy\biggr).
\end{equation}
From the properties of $u(x;t)$ in \cite{DZ2}, it is easy 
to see that $F(x;t)$, $0<t\le 1$, is indeed a distribution function.
One of the main results in \cite{BDJ} is that for fixed $x\in\R$, 
\begin{equation}\label{eq7}
\lim_{N\to\infty} 
\Prob\biggl(\chi^{(1)}_N:=\frac{l_N^{(1)}-2\sqrt{N}}{N^{1/6}}
\le x\biggr) 
=F(x;1).
\end{equation}
Convergence of the moments of $\chi^{(1)}_N$is also proved in \cite{BDJ},  
implying that 
\begin{equation}
      \lim_{N\to\infty} \frac{\Var\bigl(l^{(1)}_N\bigr)}{N^{1/3}} =
\int_{-\infty}^{\infty} t^2dF(t;1) -
\biggl( \int_{-\infty}^{\infty} tdF(t;1) \biggr)^2
=0.8132\cdots.
  \end{equation}
and
\begin{equation}
  \lim_{N\to\infty} \frac{
\Exp\bigl(l^{(1)}_N\bigr)-2\sqrt{N}}{N^{1/6}} =
\int_{-\infty}^{\infty} tdF(t;1)
=-1.7711\cdots.
\end{equation}

The above results are intimately connected to 
the Gaussian Unitary Ensemble 
(GUE) of random matrix theory.
In GUE, one considers $N\times N$ Hermitian matrices, 
and the probability density for the eigenvalues in 
an infinitesimal multi-interval about points 
$x_1,\cdots,x_N$ is given by 
\begin{equation}
   Z_N^{-1}e^{-\sum_{j=1}^N x_j^2}\prod_{1\le j<k\le N}|x_j-x_k|^2
dx_1\cdots dx_N,
\end{equation}
where $Z_N$ is the normalization constant.
Let $\lambda_{1^{\text{st}}}(M)$ be the largest eigenvalue of $M$ in GUE.
In 1994 in \cite{TW}, Tracy and Widom showed that if one scales  
$\lambda^{\text{sc}}_{1^{\text{st}}}
:=\bigl(\lambda_{1^{\text{st}}} -\sqrt{2N}\bigr)\sqrt{2}N^{1/6}$, 
then for any fixed $x\in\R$, 
\begin{equation}\label{eq11}
   \lim_{N\to\infty} 
\Prob\bigl(\lambda^{\text{sc}}_{1^{\text{st}}}\le x\bigr)
= F(x;1).
\end{equation}
In other words, properly centered and scaled, the length of the longest 
increasing subsequence for $\pi\in S_N$ 
(or $l_N^{(1)}$ under Plancherel measure on $Y_N$),
behaves statistically for large $N$ like the largest 
eigenvalue of a random GUE matrix.
Moreover, in \cite{TW} the authors also computed 
the limiting distribution of 
the $k^{\text{th}}$ largest eigenvalue 
$\lambda_{k^{\text{th}}}$ of a random GUE matrix 
for general $k\geq 2$.
Indeed, again scaling $\lambda^{\text{sc}}_{k^{\text{th}}}
:=\bigl(\lambda_{k^{\text{th}}} -\sqrt{2N}\bigr)\sqrt{2}N^{1/6}$, 
Tracy and Widom showed that 
for fixed $x\in\R$, and for $k\geq 2$,
\begin{equation}\label{eq12}
   \lim_{N\to\infty}
\Prob\bigl(\lambda^{\text{sc}}_{k^{\text{th}}}\le x\bigr)
= \lim_{N\to\infty} 
\Prob\bigl(\lambda^{\text{sc}}_{(k-1)^{\text{th}}}\le x \bigr)
+ \frac1{(k-1)!}\biggl( -\frac{\partial}{\partial t}
\biggr)^{(k-1)} \biggr|_{t=1} F(x;t).
\end{equation}

In this paper, we prove that $l^{(2)}_N$ 
behaves statistically, as $N\to\infty$,  
like the second largest eigenvalue of a random GUE matrix.
This result was conjectured in \cite{BDJ} : the obvious analogous 
result should be true for all the rows.
For the second row, the conjecture was strongly supported by 
Monte Carlo simulations of Odlyzko and Rains.
To compute the asymptotics of $l^{(2)}_N$, 
we first obtain an expression 
for the generating function $\phi_n^{(2)}(\lambda)$ 
in terms of a derivative of a Fredholm determinant.
For the case $k=1$ in \cite{BDJ}, it was already known 
(see, for example, \cite{Ge}, \cite{Ra}) 
that the generating function $\phi_n^{(1)}(\lambda)$ is a Toeplitz 
determinant :  
\begin{equation}\label{eq13}
   \phi_n^{(1)}(\lambda) = e^{-\lambda}\det\bigl(T_{n-1}\bigr), 
\qquad n\geq 0,
\end{equation}
where $T_{n-1}=\bigl((T_{n-1})_{j,k}\bigr)_{0\le j,k\le n-1}$ 
is the $n\times n$ Toeplitz matrix with respect to the weight 
$e^{2\sqrt\lambda\cos\theta}d\theta/(2\pi)$, 
\begin{equation}\label{eq-14}
   T_{n-1} = \biggl( \int_{0}^{2\pi} e^{-i(j-k)\theta}
e^{2\sqrt\lambda\cos\theta}
\frac{d\theta}{2\pi} \biggr)_{0\le j,k\le n-1}, 
\qquad n\geq 1,  
\end{equation}
and $T_{-1}$ is the $1\times 1$ matrix with entry equal to $1$.
Therefore the main part of 
the calculation in \cite{BDJ} was the computation of the 
asymptotics of $\phi_n^{(1)}(\lambda)$ as $\lambda,n\to\infty$.
In this paper, however, we focus on deriving an appropriate 
expression for the generating function $\phi_n^{(2)}(\lambda)$.
As we will explain, the asymptotics as $\lambda, n\to\infty$ 
can then be obtained in a similar manner to \cite{BDJ}.
An interesting by-product of the calculation 
of $\phi^{(2)}(\lambda)$ is a new expression for 
$\phi^{(1)}(\lambda)$ (see \eqref{eq17} below).

We state our results.
Set
\begin{equation}\label{e4}
   \varphi(z):= e^{\sqrt\lambda(z-z^{-1})}.
\end{equation}
and let $\Sigma$ be the unit circle in the complex plane 
oriented counterclockwise.
Let $\K_n$ be the integral operator acting on $L^2(\Sigma, |dw|)$,
whose kernel is defined by
\begin{equation}\label{e-16}
   \K_n(z,w) := \frac{z^{-n}w^n- 
\varphi(z)\varphi(w)^{-1}}{2\pi i(z-w)},
\qquad   \bigl(\K_n f\bigr)(z) = \int_\Sigma \K_n(z,w)f(w) dw.
\end{equation}

\begin{thm}\label{thm-1}
   For $n \geq 0$, we have
\begin{equation}\label{eq17}
   \phi^{(1)}_n(\lambda) = 2^{-n}\det(I-\K_n),
\end{equation}
and 
\begin{equation}\label{eq16}
   \phi^{(2)}_{n+1}(\lambda) = \phi^{(1)}_n(\lambda) +
\biggl( -\frac{\partial}{\partial t}\biggr) \biggr|_{t=1} \biggl[
(1+\sqrt{t})^{-n}\det(I-\sqrt{t}\K_n) \biggr].
\end{equation}
\end{thm}

\begin{thm}\label{thm-2}
   For fixed $x\in\R$, we have
\begin{equation}\label{e-21}
   \lim_{N\to\infty} 
\Prob\biggl(\chi^{(2)}_N:=\frac{l_N^{(2)}-2\sqrt{N}}{N^{1/6}}
\le x\biggr) = F^{(2)}(x),
\end{equation}
where $F^{(2)}(x)=F(x;1)+ \bigl( -\frac{\partial}{\partial t}\bigr) 
\bigr|_{t=1} F(x;t)$ is the Tracy-Widom second eigenvalue distribution 
formula given in \eqref{eq12} above.
\end{thm}

As in \cite{BDJ}, we also have convergence of the moments.
Indeed, let $\chi^{(2)}$ be a random variable 
with distribution function $F^{(2)}$.
The we have the following result.
\begin{thm}\label{thm-3}
For $m=1,2,\cdots$, 
\begin{equation}
   \lim_{N\to\infty} \Exp\bigl( \bigl(\chi^{(2)}_N\bigr)^m\bigr) 
= \Exp\bigl( \bigl(\chi^{(2)}\bigr)^m\bigr).
\end{equation}
\end{thm}
These results show that $l_N^{(2)}$ 
behaves statistically for large $N$ like the second 
largest eigenvalue of a random 
GUE matrix, under appropriate centering and scaling.
Furthermore, in view of the previous remarks, by the Robinson-Schensted 
correspondence these results show that, after appropriate centering 
and scaling, the difference between the length of the longest $2$-increasing 
subsequence and the length of the 
longest increasing subsequence of a random permutation 
of $N$ numbers, also behaves statistically for large $N$ like the second 
largest eigenvalue of a random GUE matrix.

As indicated above, the bulk of the paper is devoted to proving 
formula \eqref{eq16} for $\phi^{(2)}_{n+1}(\lambda)$.
We indicate briefly in Section \ref{sec:asymp} below 
how Theorems \ref{thm-2} and \ref{thm-3} then follow 
by the Riemann-Hilbert/steepest descent methods of \cite{BDJ} : 
further details on the computations will appear in a later 
publication.
In \cite{BDJ}, the authors express $\phi^{(1)}_n(\lambda)$ via 
\eqref{eq13} in terms of the solution of a Riemann-Hilbert problem 
(RHP) for polynomials orthogonal with respect to the weight 
$e^{2\sqrt{\lambda}\cos\theta}d\theta/(2\pi)$ on the unit circle 
(this RHP is the analog for orthogonal polynomials on the circle 
of the RHP introduced in \cite{FIK} for orthogonal polynomials on the line),  
and then apply the steepest descent method for RHP's  
introduced by Deift and Zhou in \cite{DZ1}, further developed 
in \cite{DZ2}, \cite{DVZ1}, and finally placed in a systematic form 
by Deift, Venakides and Zhou in \cite{DVZ2}, to compute the asymptotics 
as $\lambda, n\to\infty$.
A general reference for RHP's is, for example, \cite{CG}.
The calculations in \cite{BDJ} have many similarities to the 
calculations in \cite{DKMVZ}.
The methods of \cite{BDJ} apply here because the operator 
$\K_n$ in \eqref{e-16} is an example of a so-called integrable 
operator, whose resolvent can be computed in terms of a canonically 
associated RHP (see Section \ref{sec:RHP} below).
Integrable operators were introduced as a distinguished class by 
Its, Izergin, Korepin and Slavnov in \cite{IIKS}, and have since 
been applied to a broad and rapidly growing array of 
problems in pure and applied mathematics 
(see, for example, \cite{D}).

The proof of Theorem \ref{thm-1} is based, in large part, 
on manipulations of RHP's.
The proof of equality \eqref{eq17} for 
$\phi^{(1)}_n(\lambda)$ is given in Section \ref{sec:RHP} 
(see Proposition \ref{prop-1}).
In Section \ref{sec:RHP}, the RHP's $(v_Y(\cdot\thinspace;k),\Sigma)$ 
and $(v(\cdot\thinspace;k;t),\Sigma)$ are also introduced and 
their connections to the Fredholm determinant 
$\det(I-\K_N)$ are established.
In Section \ref{sec:int}, to prove \eqref{eq16} for 
$\phi^{(1)}_{n+1}(\lambda)$, we first use 
the Frobenius-Young formula 
for $d_\mu$ in \eqref{eq1} to obtain an intermediate form for  
$\phi^{(2)}_{n+1}(\lambda)$ in terms of the inverse of 
a Toeplitz matrix together with certain binomial sums 
(see Proposition \ref{lem1}).
The calculation of this intermediate form is similar to the 
derivation of \eqref{eq13} above in the Appendix in \cite{BDJ}.
The identification of this intermediate form with the right hand side of 
\eqref{eq16} is made first when $n=0$ (Section \ref{sec:0}), 
and then for general $n\geq 1$ (Section \ref{sec:positive}).
Finally, in Section \ref{sec:asymp}, we indicate how to prove 
Theorems \ref{thm-2} and \ref{thm-3} following the methods in 
\cite{BDJ}. 
In the Appendix, we discuss the spectral properties of the 
operator $\K_n$ in \eqref{e-16}, and the 
(unique) solvability of the RHP \eqref{e+5} below.

The proof that we give for the basic formula \eqref{eq16} 
appears rather ad hoc.
At the end of the paper, in Section \ref{sec:asymp}, 
we provide a motivation for our calculations.\\

\noindent {\bf Acknowledgments.} 
The authors would like to thank Andrew Odlyzko and Eric Rains
for their Monte Carlo simulations on the statistics of the
second row. The authors also greatly appreciate many useful conversations
and communications with Persi Diaconis, Andrew Odlyzko, Eric Rains,
Craig Tracy and Harold Widom. The work of the authors was supported
in part by a Sloan Doctoral Dissertation Fellowship [J.B.],
by NSF grant \#DMS-9500867 and the Institute of Advanced Study [P.D.], 
and by the Swedish National Research Council (NFR) [K.J.].


\section{Riemann-Hilbert problems (RHP's)}\label{sec:RHP}
As above, let $\Sigma$ denote the unit 
circle in the complex plane, oriented counterclockwise.
Set 
\begin{equation}\label{e+1}
   \psi(z) := e^{\sqrt{\lambda}(z+z^{-1})}.
\end{equation}
In this paper, we use the following two (matrix) RHP's.
First, for $k\geq 0$, let $Y(z;k)$, a $2\times 2$ matrix-valued 
function of $z$, 
be the solution of the RHP $(v_Y(\cdot\thinspace;k),\Sigma)$,
\begin{equation}\label{eq-19}
 \begin{cases}
     Y(z;k) \quad\text{is analytic in}\quad z\in\C\setminus\Sigma,\\
     Y_+(z;k)=Y_-(z;k) \begin{pmatrix}
1&\frac{1}{z^{k}}\psi(z)\\0&1 \end{pmatrix},
\quad\text{on}\quad z\in\Sigma,\\
     Y(z;k) \biggl(\begin{smallmatrix} z^{-k}&0\\0&z^k
\end{smallmatrix} \biggr)
=I+O(\frac1{z}) \quad\text{as}\quad z\to\infty.
 \end{cases}
\end{equation}
The notation $Y_+(z;k)$ (resp., $Y_-$) denotes the limiting value 
$\lim_{z'\to z}Y(z';k)$ with $|z'|<1$ (resp., $|z'|>1$).
Note that $k$ plays the role of an external parameter in 
\eqref{eq-19} ; in particular, the term 
$O(\frac1{z})$ does not imply a uniform bound in $k$.
For a general RHP, the existence and uniqueness of the solution 
is, of course, not clear a priori. 
But for the case at hand, we can simply write down the 
(unique) solution explicitly in terms 
of orthogonal polynomials, as follows.

Let $\pi_n(z)=z^n+\sum_{p=0}^{n-1}\eta^n_{p}z^{p}$ 
denote the $n$-th monic orthogonal polynomial with 
respect to the measure $\psi(z)dz/(2\pi iz)$ 
on the unit circle, and introduce the polynomial 
$\pi_n^{*}(z):=z^n\overline{\pi_n}(1/z)
=z^n(z^{-n}+\sum_{p=0}^{n-1}\overline{\eta^n_{p}}z^{-p}$ (see \cite{Sz}). 
For the measure at hand, $\psi(z)dz/(2\pi iz)
=e^{2\sqrt\lambda\cos\theta} d\theta/(2\pi)$ , 
all the coefficients of $\pi_n(z)$ are real 
and $\pi_n^{*}(z)=z^n\pi_n(1/z)$. 
The solution of the RHP $(v_Y(\cdot\thinspace;k),\Sigma)$ is given by 
the following formula (see Lemma 4.1 in \cite{BDJ} ;
here we change the notation $f(z)$ to $\psi(z)$ and use the 
monic orthogonal polynomial $\pi_k(z)$ instead of 
the normalized orthogonal polynomials $p_k(z)$),
\begin{equation}\label{e+3}
    Y(z;k)=\begin{pmatrix}
    \pi_{k}(z)&
\int_{\Sigma}\frac{\pi_{k}(s)}{s-z}\frac{\psi(s)ds}{2\pi i s^{k}}\\
    -\kappa^2_{k-1}\pi^{*}_{k-1}(z)& -\kappa^2_{k-1}
\int_{\Sigma}\frac{\pi^{*}_{k-1}(s)}{s-z}\frac{\psi(s)ds}{2\pi i s^{k}}
       \end{pmatrix}, \qquad k\geq 1.
\end{equation}
For $k=0$, the solution is given by \eqref{eq-73} below.

Let $\varphi(z)$ be defined as in \eqref{e4}.
For $0<t\le 1$ and for $k\geq 0$,  
the second RHP 
$(v(\cdot\thinspace;k;t),\Sigma)$ is to find $m(z;k;t)$ satisfying 
\begin{equation}\label{e+5}
  \begin{cases}
   &m(z;k;t) \quad\text{is analytic in} \quad z\in\C\setminus\Sigma,\\
   &m_+(z;k;t)=m_-(z;k;t)\begin{pmatrix}
1-t&-\sqrt{t}z^{-k}\varphi(z)^{-1}\\
\sqrt{t}z^{k}\varphi(z)&1 \end{pmatrix}
\quad \text{on $z\in\Sigma$},\\
   &m(z;k;t)\to I \quad\text{as} \quad z\to\infty.
  \end{cases}
\end{equation}
In this RHP, there are two external parameters $k$ and $t$.
For $t=1$, this RHP is equivalent to the RHP 
$(v_Y(\cdot\thinspace;k),\Sigma)$, in the sense that a solution of one RHP 
implies a solution of the other RHP.
Indeed one can easily check that
\begin{equation}\label{e+6}
   m(z;k;1)=
   \begin{cases}
\biggl( \begin{smallmatrix} 0&-1\\1&0 \end{smallmatrix} \biggr)Y(z;k)
\biggl( \begin{smallmatrix} e^{\sqrt{\lambda}z}&0\\
0&e^{-\sqrt{\lambda}z} \end{smallmatrix} \biggr),  &|z|<1,\\
\biggl( \begin{smallmatrix} 0&-1\\1&0 \end{smallmatrix} \biggr) Y(z;k)
\biggl( \begin{smallmatrix} z^{-k}e^{\sqrt{\lambda}z^{-1}}&0\\
0&z^{k}e^{-\sqrt{\lambda}z^{-1}} \end{smallmatrix} \biggr)
\biggl( \begin{smallmatrix} 0&1\\-1&0 \end{smallmatrix} \biggr), 
\quad &|z|>1,
 \end{cases}
\end{equation}
solves \eqref{e+5} if and only if $Y(z;k)$ solves \eqref{eq-19}. 
The (unique) solvability of the RHP \eqref{e+5} for $0<t<1$ 
is discussed in the Appendix.
As noted earlier, 
the motivation for the introduction of \eqref{e+5} for $t<1$ is given 
in Section \ref{sec:asymp}  below.

The solution $m(z;k;t)$ of the RHP $(v(z;k;t),\Sigma)$ is 
related to the Fredholm determinant of the integral operator 
$\K_n$ in \eqref{e-16} acting on $L^2(\Sigma, |dw|)$. 
For each $n$, it is easy to verify that $\K_n$ is trace class 
and hence $\det(I-\K_n)$ exist.
The proof of the following Lemma is modeled on 
Proposition 6.13 in \cite{DIZ}.

\begin{lem}\label{lem-1}
   Let $m(z;k;t)$ be the solution of the RHP $(v(z;k;t),\Sigma)$ 
given in \eqref{e+5}, and let $\K_k$ be defined as above.
If we denote the 11-component of $m(z;k;t)$ by $m_{11}(z;k;t)$, 
we have at $z=0$, 
\begin{equation}
   m_{11}(0;k;t)=(1+\sqrt{t})
\frac{\det(I-\sqrt{t}\K_{k-1})}{\det(I-\sqrt{t}\K_k)}, 
\qquad k\geq 0.
\end{equation}
\end{lem}

\begin{proof}
The operator $\K_n$ has norm less than or equal to $1$ and $1$ is not 
an eigenvalue of $\K_n$ (see Appendix).
So $\det(I-\sqrt{t}\K_k)$ never vanishes for any $k$, $0<t\le 1$.

First note that
\begin{equation}
   \K_{k-1}(z,w) = \K_k(z,w) + \frac{1}{2\pi i}z^{-k}w^{k-1}
=: \K_k(z,w)+E_k(z,w).
\end{equation}
Then we have
\begin{equation}
   \frac{\det(I-\sqrt{t}\K_{k-1})}{\det(I-\sqrt{t}\K_k)}
=\det\biggl(I-\frac1{I-\sqrt{t}\K_k}\sqrt{t}E_k\biggr).
\end{equation}
Since $E_k$ is a rank 1 operator,
\begin{equation}\label{e+15}
   \frac{\det(I-\sqrt{t}\K_{k-1})}{\det(I-\sqrt{t}\K_k)}
=1-\tr\biggl(\frac1{I-\sqrt{t}\K_k}\sqrt{t}E_k\biggr)
=1-\frac{\sqrt{t}}{2\pi i}\int_\Sigma
\biggl(\frac1{I-\sqrt{t}\K_k}f_1\biggr)(z)\cdot z^{k-1} dz,
\end{equation}
where $f_1(z)=z^{-k}$ as in \eqref{e+10} below.

On the other hand, we define $M(z,k;t)$ as follows,
\begin{equation}\label{e+7}
   M(z;k;t):=
   \begin{cases}
m(z;k;t) \biggl( \begin{smallmatrix} 1/(1+\sqrt{t})&0\\0&1+\sqrt{t}
\end{smallmatrix} \biggr),
\qquad &|z|<1\\
m(z;k;t),  &|z|>1.
 \end{cases}
\end{equation}
Then it is easy to check that $M(z;k;t)$ solves 
a new RHP $(V(\cdot\thinspace ;k;t),\Sigma)$,
\begin{equation}\label{e+7.5}
  \begin{cases}
       &M(z;k;t) \quad\text{is analytic in} \quad \C\setminus\Sigma,\\
   &M_+(z;k;t)=M_-(z;k;t)\begin{pmatrix}
1-\sqrt{t}&-\sqrt{t}(1+\sqrt{t})z^{-k}\varphi(z)^{-1}\\
\frac{\sqrt{t}}{1+\sqrt{t}}z^{k}\varphi(z)&1+\sqrt{t} \end{pmatrix}
\quad \text{on $z\in\Sigma$},\\
   &M(z;k;t)\to I \quad\text{as} \quad z\to\infty.
  \end{cases}
\end{equation}
Note that the jump matrix $V(z;k;t)$ can be written in the form 
\begin{equation}\label{eq25}
   V(z;k;t)=I-2\pi i \sqrt{t}f(z;k;t)\bigl(g(z;k;t)\bigr)^T,
\end{equation}
where $f, g$ are column vectors defined by
\begin{equation}\label{e+10}
   f(z;k;t)= (f_1,f_2)^T=\bigl(z^{-k}, 
-\frac{1}{1+\sqrt{t}}\varphi(z) \bigr)^T, \quad
g(z;k;t)= (g_1,g_2)^T=\frac{1}{2\pi i}
\bigl( z^k, (1+\sqrt{t})\varphi(z)^{-1} \bigr)^T.
\end{equation}
From the general theory of RHP's, 
the solution $M(z;k;t)$ of the RHP $(V(z;k;t),\Sigma)$ satisfies 
\begin{equation}
 \begin{split}
   M(z;k;t)&=I+\frac1{2\pi i} \int_\Sigma 
\frac{M_+(s;k;t)\bigl(I-V(s;k;t)^{-1}\bigr)}{s-z} ds\\
&=I-\frac{\sqrt{t}}{2\pi i} \int_\Sigma 
M_+(s;k;t)\biggl( \begin{smallmatrix} 1&(1+\sqrt{t})s^{-k}\varphi(s)^{-1}\\
-\frac1{1+\sqrt{t}}s^k\varphi(s)&-1 \end{smallmatrix}\biggr) \frac{ds}{s-z}.
 \end{split}
\end{equation}
Recalling the definition of $f(z;k;t)$, we have 
\begin{equation}\label{e+17}
   M_{11}(0;k;t)=1-\frac{\sqrt{t}}{2\pi i} \int_\Sigma 
\bigl( M_+(z;k;t)f(z;k;t) \bigr)_{1} 
\cdot z^{k-1}dz.
\end{equation}
Now from the theory of integrable operators (see for example \cite{D}),  
the integral operator $S_k$ acting on $L^2(\Sigma, |dw|)$ 
with kernel  
\begin{equation}
   S_k(z,w)  := \frac{(f(z;k;t))^Tg(w;k;t)}{z-w},
\end{equation}
satisfies the relation 
\begin{equation}\label{e+18}
   \biggl( \frac1{I-\sqrt{t}S_k}S_k \biggr)(z,w)
=\frac{\bigl(F(z;k;t)\bigr)^TG(w;k;t)}{z-w}, 
\end{equation}
where $F, G$ are column vectors given by 
\begin{equation}\label{e+19}
  \begin{split}
&F(z;k;t)=\biggl(\frac1{I-\sqrt{t}S_k}f\biggr)(z)=M_+(z;k;t)f(z;k;t),\quad \\
&G(w;k;t)=\biggl(\frac1{I-\sqrt{t}S_k}g\biggr)(w)
=\bigl(M_+(w;k;t)^T\bigr)^{-1}g(w;k;t).
  \end{split}
\end{equation}
But from the definitions of $f$ and $g$ in \eqref{e+10},
\begin{equation}\label{eq40}
   S_k(z,w)  = \K_k(z,w). 
\end{equation}
Therefore \eqref{e+15}, \eqref{e+17} and the first part of \eqref{e+19} 
give us that 
\begin{equation}
  M_{11}(0;k;t)=\frac{\det(I-\sqrt{t}\K_{k-1})}{\det(I-\sqrt{t}\K_k)}.
\end{equation}
The relation \eqref{e+7} completes the proof.
\end{proof}

Now for fixed $t$, $\lambda$, we compute the asymptotics of 
$\det(I-\sqrt{t}\K_p)$ as $p\to\infty$. 
The calculation below is similar to that of \cite{D} where a new proof 
of Szeg\"o's strong limit theorem is given.
\begin{lem}\label{lem-2}
   For fixed $0<t\le 1$ and $\lambda>0$, we have 
\begin{equation}
   \lim_{p\to\infty} (1+\sqrt{t})^{-p}\det(I-\sqrt{t}\K_p) =1.
\end{equation}
\end{lem}

\begin{proof}
   First note that 
\begin{equation}\label{e+22}
   \log\det(I-\sqrt{t}\K_p) = -\int_0^{\sqrt{t}} 
\tr\biggl( \frac1{I-s\K_p}\K_p \biggr) ds.
\end{equation}
From \eqref{e+18} and \eqref{eq40}, the right hand side 
of \eqref{e+22} is given by 
\begin{equation}\label{e+23}
   -\int_0^{\sqrt{t}} \tr\biggl( \frac{\bigl(F(z;p;s^2)\bigr)^TG(w;p;s^2)}{z-w} 
\biggr) ds
= -\int_0^{\sqrt{t}} ds \int_{\Sigma} \bigl(F'(z;p;s^2)\bigr)^TG(z;p;s^2) dz
\end{equation}
where the prime $'$ indicates differentiation with respect to $z$.
From \eqref{e+19}, in order to compute the asymptotics of $\det(I-\sqrt{t}\K_p)$, 
we need the asymptotics of $M_+(z;p;s^2)$ as $p\to\infty$ 
uniformly in $0<s\le \sqrt{t}$.

Note that the jump matrix $V(z;p;s^2)$ has the following factorization 
\begin{equation}
   V(z;p;s^2) = \begin{pmatrix} 1&-sz^{-p}\varphi^{-1}\\0&1\end{pmatrix} 
\begin{pmatrix} \frac1{1+s}&0\\0&1+s \end{pmatrix} 
\begin{pmatrix} 1&0\\ \frac{s}{(1+s)^2}z^p\varphi &1\end{pmatrix}.
\end{equation}
Let $0< \rho <1$ be any number.
Define 
\begin{equation}
  \begin{cases}
     \tilde{M}(z) := M(z;p;s^2), & |z|<\rho, \\
     \tilde{M}(z) := M(z;p;s^2) \begin{pmatrix} 
1&0\\ \frac{s}{(1+s)^2}z^p\varphi &1\end{pmatrix}^{-1}, &\rho <|z|<1, \\
     \tilde{M}(z) := M(z;p;s^2) \begin{pmatrix} 
1&-sz^{-p}\varphi^{-1}\\0&1\end{pmatrix}, &1<|z|<\rho^{-1}, \\
     \tilde{M}(z) := M(z;p;s^2), & |z|>\rho^{-1},
  \end{cases}
\end{equation}
and set 
\begin{equation}
  \begin{cases}
     \tilde{V}(z)= \begin{pmatrix}
1&0\\ \frac{s}{(1+s)^2}z^p\varphi &1\end{pmatrix} , &|z| = \rho , \\
     \tilde{V}(z)= \begin{pmatrix}
\frac1{1+s}&0\\0&1+s \end{pmatrix} , &|z|=1, \\
     \tilde{V}(z)= \begin{pmatrix}
1&-sz^{-p}\varphi^{-1}\\0&1\end{pmatrix} , &|z|=\rho^{-1}.
  \end{cases}
\end{equation}
Then $\tilde{M}$ solves the new RHP $(\tilde{V}, \tilde{\Sigma})$ 
where $\tilde{\Sigma}:=\{ |z|=\rho \}\cup 
\{ |z|=1 \}\cup\{ |z|=\rho^{-1} \}$, 
oriented counterclockwise on each of the three circles,
\begin{equation}
  \begin{cases}
      \tilde{M}(z) \quad \text{is analytic in 
$\C\setminus \tilde{\Sigma}$,} \\
       \tilde{M}_+(z)= \tilde{M}_-(z) \tilde{V}(z), 
\quad\text{on $z\in \tilde{\Sigma}$,} \\
       \tilde{M}(z) \to I, \qquad \text{as $z\to\infty$}.
  \end{cases}
\end{equation}
Here $\tilde{M}_+(z)= \lim_{z'\to z, |z'|<|z|} \tilde{M}(z')$ 
and $\tilde{M}_-(z)= \lim_{z'\to z, |z'|>|z|} \tilde{M}(z')$ 
for $z$ on each of the three circles.
Now as $p\to\infty$, we see that $\tilde{V}(z)\to I$ for 
$z\in\{ |z|=\rho \}\cup\{ |z|=\rho^{-1} \}$.
Moreover the convergence is exponential at a rate which is 
uniform in $z$ and $s$.
Therefore it follows that $\tilde{M}\to \tilde{M}^\infty$ where 
$\tilde{M}^\infty$ is the solution of the RHP $(\tilde{V}^\infty, \Sigma)$, 
$\tilde{V}^\infty := \biggl( \begin{smallmatrix}
\frac1{1+s}&0\\0&1+s \end{smallmatrix} \biggr)$.
This RHP has the explicit solution 
\begin{equation}
  \begin{cases}
     \tilde{M}^\infty = \biggl( \begin{smallmatrix}
\frac1{1+s}&0\\0&1+s \end{smallmatrix} \biggr), &|z|<1, \\
     \tilde{M}^\infty = I, &|z|>1.
  \end{cases}
\end{equation}
Therefore as $p\to\infty$, we have the following 
asymptotics for $z\in \Sigma$, 
\begin{equation}
  M_+(z;p;s^2)=\tilde{M}_+(z) \begin{pmatrix}
1&0\\ \frac{s}{(1+s)^2}z^p\varphi &1\end{pmatrix}
\sim \tilde{M}_+^\infty(z) \begin{pmatrix}
1&0\\ \frac{s}{(1+s)^2}z^p\varphi &1\end{pmatrix}
= \begin{pmatrix} \frac1{1+s}&0\\ \frac{s}{1+s}z^p\varphi &1+s
\end{pmatrix} ,
\end{equation}
where the error is exponentially small, uniformly for $z$ and $s$.

Now from \eqref{e+19}, the above asymptotic result yields 
for $z\in\Sigma$, 
\begin{equation}
 \begin{split}
   &F(z;p;s^2)= M_+(z;p;s^2)f(z;p;s^2) \sim 
\biggl( \frac{z^{-p}}{1+s}, \frac{-\varphi}{1+s} \biggr)^T,\\
   &F'(z;p;s^2) \sim \biggl( \frac{-pz^{-p-1}}{1+s}, 
\frac{-\sqrt\lambda(1+z^{-2})\varphi}{1+s} \biggr)^T,\\
   &G(z;p;s^2)= \bigl(M_+(z;p;s^2)^T\bigr)^{-1}g(z;p;s^2) \sim 
\frac1{2\pi i}\bigl( z^p, \varphi^{-1} \bigr)^T.
 \end{split}
\end{equation}
Inserting these asymptotics, which have exponential error terms,  
into \eqref{e+23}, equation \eqref{e+22} becomes 
\begin{equation}
 \begin{split}
   \log\det(I-\sqrt{t}\K_p) &= 
-\int_0^{\sqrt{t}} \frac{ds}{2\pi i(1+s)} \int_{\Sigma}  
\bigl[ -pz^{-1}-\sqrt\lambda(1+z^{-2}) \bigr]  dz
+ O\bigl(e^{-cp}\bigr) \\
&=p\log(1+\sqrt{t}) + O\bigl(e^{-cp}\bigr), 
 \end{split}
\end{equation}
where the exponential error term is uniform for $0<t\le 1$.
\end{proof}

Combining Lemma \ref{lem-1} and Lemma \ref{lem-2}, we summarize the 
results in this section as follows.
\begin{prop}\label{prop-1}
   Let $m(z;k;t)$ be the solution of the RHP 
$(v(\cdot\thinspace ;t),\Sigma)$ given
in \eqref{e+5}, and let $\K_n$ be defined as in \eqref{e-16}.
Denoting the $11$-component of $m(z;k;t)$ 
by $m_{11}(z;k;t)$, we have
\begin{equation}\label{eq-49}
   \prod_{k=n}^{\infty} m_{11}(0;k+1;t) = (1+\sqrt{t})^{-n}
\det(I-\sqrt{t}\K_n), \qquad n\geq 0.
\end{equation}
In particular, we have,  
for the generating function of $l_N^{(1)}$,
\begin{equation}\label{eq-50}
   \phi^{(1)}_n(\lambda) 
= 2^{-n} \det(I-\K_n), \qquad n\geq 0.
\end{equation}
Also for $n\geq 0$, the right hand side of \eqref{eq16} is given by 
\begin{equation}\label{eq49}
   \phi^{(1)}_n(\lambda) +
\biggl( -\frac{\partial}{\partial t}\biggr) \biggr|_{t=1} \biggl[
(1+\sqrt{t})^{-n}\det(I-\sqrt{t}\K_n) \biggr]=
   \biggl[ 1-\sum_{k=n}^{\infty} \frac{\dot{m}_{11}(0;k+1;1)}{m_{11}(0;k+1;1)} 
\biggr] \phi^{(1)}_n(\lambda),
\end{equation}
where the dot $\cdot$ indicates differentiation with respect to $t$.
\end{prop}

\begin{proof}
  Equation \eqref{eq-49} follows from Lemma \ref{lem-1} and Lemma \ref{lem-2}.

  From (1.25) and (1.27) in \cite{BDJ}, we have 
\begin{equation}
   \phi^{(1)}_n(\lambda) = \prod_{k=n}^{\infty} \bigl(-Y_{21}(0;k+1) \bigr),
\end{equation}
where $Y(z;k)$ solves the RHP \eqref{eq-19}.
From \eqref{e+6}, $m_{11}(0;k;1)=-Y_{21}(0;k)$.
Hence we have 
\begin{equation}
   \phi^{(1)}_n(\lambda) = \prod_{k=n}^{\infty} m_{11}(0;k+1) 
= 2^{-n} \det(I-\K_n).
\end{equation}

Equation \eqref{eq49} is obtained by taking derivative of \eqref{eq-49} 
and using \eqref{eq-50}.
(The differentiability of the infinite product in \eqref{eq-49} follows from 
the uniform error estimate in the proof of Lemma \ref{lem-2}.)
\end{proof}


\section{Intermediate form of the generating function}\label{sec:int}
For $N,n \geq 0$, let $q^{(2)}_{n,N}$ and $\phi^{(2)}_n(\lambda)$ be defined 
as in \eqref{e-1} and \eqref{e-2} with $k=2$.
Define 
\begin{equation}\label{e+35}
  \begin{split}
   &a_0(s) := 
\sum_{m=1}^{\infty} \frac{\lambda^{m}}{(m+s)^2((m-1)!)^2}, \qquad s\geq 1,\\
   &a_0(0) := \sum_{m=0}^{\infty} \frac{\lambda^{m}}{(m!)^2}, \\
   &a_j(s) := \sum_{m=j}^{\infty} \frac{\lambda^{m-j/2}}{(m+s)(m-1)!(m-j)!},
\qquad s\geq 0, j\geq 1\\
   &b_n(s) := \bigl(a_1(s), \cdots, a_{n}(s)\bigr)^T, \qquad s\geq 0.
  \end{split}
\end{equation}
Let $\T_{n-1}$ be the $n\times n$ Toeplitz matrix given in \eqref{eq-14},  
\begin{equation}\label{eq58}
   \T_{n-1} = \biggl( \int_{0}^{2\pi} e^{-i(j-k)\theta}
e^{2\sqrt\lambda\cos\theta} 
\frac{d\theta}{2\pi} \biggr)_{0\le j,k\le n-1}
\end{equation}

\begin{prop}\label{lem1}
We have the following expressions for $\phi^{(2)}_{n+1}(\lambda)$, 
\begin{eqnarray}\label{eq52}
   \phi^{(2)}_{n+1}(\lambda)& =& \biggl[ \sum_{s=0}^{\infty} 
\frac{\lambda^{s}}{(s!)^2} \biggl( a_{0}(s)- 
\bigl( T^{-1}_{n-1}b_{n}(s),b_{n}(s)\bigr)\biggr)\biggr] 
\phi^{(1)}_n(\lambda), \qquad n\geq 1,\\
\label{eq52+} \phi^{(2)}_{1}(\lambda)& =& \biggl[ \sum_{s=0}^{\infty}
\frac{\lambda^{s}}{(s!)^2} a_{0}(s) \biggr] \phi^{(1)}_0(\lambda).
\end{eqnarray}
\end{prop}

\begin{proof}
First note that the statistics of the second row of $\mu\in Y_N$ is the same 
as the second column of $\mu\in Y_N$. 
This follows immediately by considering the fact that if 
$\pi$ has a $k$-increasing subsequence, then $\pi'(j):=N+1-\pi(j)$ has a 
$k$-decreasing subsequence, together with the combinatorial 
results of Greene referred to in the Introduction relating the rows 
and columns of a Young tableaux to $k$-increasing and 
$k$-decreasing subsequences.
Alternatively, if $\mu\in Y_N$ and if $\mu'$ is its transpose, 
then the hook formula implies that $d_{\mu}=d_{\mu'}$, which in turn yields 
that the $k$-th row and the $k$-th column have the same statistics 
under Plancherel measure.

Let $n\geq 1$. 
For a partition $\mu\vdash N$, $r_1(\mu)$ and $r_2(\mu)$ are 
defined to be the lengths of the first and 
the second columns, respectively. 
Observe that $\mu = (\mu_1, \mu_2, \cdots, \mu_{r_2(\mu)}, 
1,1,\cdots,1)$ and $\mu_1\geq \mu_2\geq \cdots\geq \mu_{r_2(\mu)} \geq 2$.
From the remark above, we have 
\begin{equation}\label{e-3}
   q^{(2)}_{n,N} = 
\sum_{\substack{\mu\vdash N \\ l^{(2)}_N(\mu)\le n}} \frac{d_\mu^2}{N!} 
= \sum_{\substack{\mu\vdash N \\ r_2(\mu)\le n}} \frac{d_\mu^2}{N!} 
= \sum_{\substack{\mu\vdash N \\ r_1(\mu)\le n}} \frac{d_\mu^2}{N!}
+ \sum_{r=n+1}^{N}\sum_{\substack{\mu\vdash N \\ r_1(\mu)=r \\ 
\mu_{n+1}=\cdots =\mu_r=1}} \frac{d_\mu^2}{N!}.
\end{equation}
Set $h_j=\mu_j+r_1(\mu)-j$.
Then we have the Frobenius-Young determinant formula (see, for example, 
\cite{Kn} 5.1.4 (34))
\begin{equation}\label{e-4}
   d_\mu = N! \prod_{1\le i<j\le r_1(\mu)}(h_i-h_j) 
\prod_{j=1}^{r_1(\mu)}\frac1{h_j!},
\end{equation}
where $\prod_{1\le i<j\le r_1(\mu)}(h_i-h_j) :=1$ for $r_1(\mu)=1$.
For $\mu\vdash N$ satisfying $r_1(\mu)=r \geq n+1$ 
and $\mu_{n+1}=\cdots =\mu_r=1$, 
the above expression becomes, after some elementary calculations,  
\begin{equation}\label{e-5}
   d_\mu = \frac{N!}{(r-n)!} \prod_{i=1}^{n}\frac{(h_i-1)!}{(h_i-r+n-1)!}
\prod_{1\le i<j\le n}(h_i-h_j) \prod_{j=1}^{n}\frac1{h_j!}.
\end{equation}
Then \eqref{e-3} may be re-expressed as 
\begin{equation}\label{e-6}
   q^{(2)}_{n,N} = q^{(1)}_{n,N}
\quad  + N!\sum_{r=n+1}^{N} \frac1{((r-n)!)^2} 
\sum_{ (*)} \prod_{1\le i<j\le n}(h_i-h_j)^2 
\prod_{j=1}^{n}\frac1{(h_j)^2((h_j-r+n-1)!)^2}
\end{equation}
where $(*)$ means that we sum over all integers 
$h_1>h_2>\cdots >h_n\geq 1+r-n$ 
such that $\sum_{j=1}^{n}h_j=N+(r-\frac12 n)(n-1)$.
Of course, the sum in \eqref{e-6} must be replaced 
by $0$ if $n\geq N$. 
Also, as above, $\prod_{1\le i<j\le n}(h_i-h_j)^2 :=1$ if $n=1$.
Hence \eqref{e-2} becomes 
\begin{equation}\label{e-7}
  \begin{split}
   &\phi^{(2)}_n(\lambda) = \phi^{(1)}_n(\lambda) \\
&\quad + e^{-\lambda}\sum_{r=n+1}^{\infty} \frac1{((r-n)!)^2}
\sum_{N\geq r}\sum_{ (*) } \lambda^N 
\prod_{1\le i<j\le n}(h_i-h_j)^2 
\prod_{j=1}^{n}\frac1{(h_j)^2(h_j-r+n-1)!)^2}.
  \end{split}
\end{equation}
Now rewriting $\lambda^N = \lambda^{-(r-\frac12 n)(n-1)+\sum_{j=1}^nh_j}$ 
and changing $\sum_{N\geq r}\sum_{ (*) } $ into the sum over all integers 
$h_1>h_2>\cdots >h_n\geq 1+r-n$ without any sum restriction, 
the above expression becomes  
\begin{equation}\label{e-8}
 \begin{split}
  &\phi^{(2)}_n(\lambda) = \phi^{(1)}_n(\lambda)\\
&\quad + e^{-\lambda}\sum_{r=n+1}^{\infty} 
\frac{\lambda^{-(r-\frac12 n)(n-1)}}{((r-n)!)^2}
\sum_{h_1>\cdots >h_n \geq 1+r-n} 
\biggl[ \prod_{1\le i<j\le n}(h_i-h_j)^2 
\prod_{j=1}^{n}\frac{\lambda^{h_j}}{(h_j)^2(h_j-r+n-1)!)^2} 
\biggr].
 \end{split}
\end{equation}
Changing the summation index to $s:=r-n$, setting $l_j:=h_j-s$, 
using the symmetry of the summand under permutation of the $h_j$'s, 
and noting that the strict ordering of the $h_j$'s is automatically 
enforced, we have 
\begin{equation}\label{e-9}
  \phi^{(2)}_n(\lambda) = \phi^{(1)}_n(\lambda)
 + e^{-\lambda}\lambda^{-n(n-1)/2} \sum_{s=1}^{\infty} 
\frac{\lambda^{s}}{(s!)^2} H(s)
\end{equation}
where 
\begin{equation}\label{e-10}
   H(s) = \frac1{n!}\sum_{l_1,\cdots,l_n\geq 1} 
\biggl[ \prod_{1\le i<j\le n}(l_i-l_j)^2
\prod_{j=1}^{n}\frac{\lambda^{l_j}}{(l_j+s)^2((l_j-1)!)^2}
\biggr], \qquad s\geq 1.
\end{equation}
Now observe (\cite{Sz}, Chapter 2) that $H(s)$ is a Hankel 
determinant with respect to the discrete measure 
\begin{equation}\label{e-11}
   \nu_s({m})=\frac{\lambda^{m}}{(m+s)^2((m-1)!)^2}, \qquad 
m\in\{1,2,\cdots\}, \quad s\geq 1,
\end{equation}
\begin{equation}\label{e-12}
   H(s)= \det \biggl( \sum_{m=1}^{\infty} 
m^{j+k}\frac{\lambda^{m}}{(m+s)^2((m-1)!)^2}\biggr)_{0\le j,k\le n-1}, 
\quad s\geq 1.
\end{equation}
Noting that 
\begin{equation}\label{e-12.0}
  \int_0^{2\pi} e^{-i(j-k)\theta} 
e^{2\sqrt{\lambda}\cos\theta} \frac{d\theta}{2\pi} 
= \sum_{m=\max(j,k)}^\infty \frac{\lambda^{m-(j+k)/2}}{(m-j)!(m-k)!},
\end{equation}
the identity \eqref{eq13} implies 
\begin{equation}\label{e-12.1}
 \begin{split}
   \phi^{(1)}_n(\lambda) &= e^{-\lambda}\lambda^{-n(n-1)/2}
\det \biggl( \sum_{m=0}^{\infty}
m^{j+k}\frac{\lambda^{m}}{(m!)^2}\biggr)_{0\le j,k\le n-1}\\
&= e^{-\lambda}\lambda^{-n(n-1)/2} H(0),
 \end{split}
\end{equation}
where $H(0)$ is the Hankel determinant with respect to the discrete measure 
\begin{equation}
   \nu_0({m})=\frac{\lambda^{m}}{(m!)^2}, \qquad
m\in\{0,1,2,\cdots\}.
\end{equation}
Therefore we have 
\begin{equation}\label{e-12.2}
   \phi^{(2)}_n(\lambda) = e^{-\lambda}\lambda^{-n(n-1)/2} \sum_{s=0}^{\infty}
\frac{\lambda^{s}}{(s!)^2} H(s).
\end{equation}
Using elementary row and column operations, we re-express $H(s)$ as 
\begin{equation}\label{e-13}
   H(s)=\det \bigl( h(s)_{jk} \bigr)_{0\le j,k\le n-1},
\end{equation}
where $h(s)_{jk}$ are defined by 
\begin{equation}\label{e-14}
 \begin{split}
    & h(s)_{00}=\sum_{m=1}^\infty \frac{\lambda^{m}}{(m+s)^2((m-1)!)^2}, 
\qquad\qquad s\geq 1,\\
    & h(0)_{00}=\sum_{m=0}^\infty \frac{\lambda^{m}}{((m!)^2},\\
    & h(s)_{0k}= h(s)_{k0}= 
\sum_{m=k}^\infty \frac{\lambda^{m}}{(m+s)(m-1)!(m-k)!}, 
\qquad k\geq 1, \quad s\geq 0,\\
    & h(s)_{jk}=\sum_{m=\max(j,k)}^\infty \frac{\lambda^{m}}{(m-j)!(m-k)!}, 
\qquad\qquad j,k \geq 1, \quad s\geq 0. 
 \end{split}
\end{equation}
We obtain 
\begin{equation}\label{e-15}
   \phi^{(2)}_n(\lambda) =
e^{-\lambda} \sum_{s=0}^{\infty} \frac{\lambda^{s}}{(s!)^2} 
\det \bigl( \lambda^{-(j+k)/2}h(s)_{jk} \bigr)_{0\le j,k\le n-1}.
\end{equation}
Setting $n=1$ in \eqref{e-15}, we immediately obtain \eqref{eq52+}.
For $n\geq 2$, from \eqref{e-12.0}, 
$\bigl( \lambda^{-(j+k)/2}h(s)_{jk} \bigr)_{1\le j,k\le n-1}$ 
is precisely the Toeplitz matrix $T_{n-2}$ in \eqref{eq58}.
Thus for each $s$ in \eqref{e-15}, 
$\bigl( \lambda^{-(j+k)/2}h(s)_{jk} \bigr)_{0\le j,k\le n-1}$
is a rank 2 extension of $T_{n-2}$, and hence can be evaluated 
in the standard way.
Indeed expanding the determinant along the first row, and then 
expanding each of the determinants obtained along the first column, 
and using \eqref{eq13}, 
$e^{-\lambda}\det(\T_{n-2})=\phi^{(1)}_{n-1}(\lambda)$, 
we obtain \eqref{eq52} directly.
\end{proof}

\noindent {\bf Remarks.} 

\noindent  (1)  Instead of expanding each 
determinant in \eqref{e-15} along rows and columns 
as above, we can appeal directly to the formulae of Weinstein and 
Aronszajn for the determinant of a finite rank extension of 
a given operator (see, for example, \cite{Ka} Chapter 4.6).

\noindent  (2) As we will see in Lemma \ref{lem3} below, 
$T^{-1}_{n-1}$, and hence $\phi^{(2)}_{n+1}(\lambda)$, involves 
full knowledge of all the monic polynomials $\pi_k(z)$ 
of degree $k\le n-1$. 
This is in contrast to $\phi^{(1)}_n(\lambda)$, 
which involves \emph{only} the leading coefficients 
$\kappa_m$ of the normalized 
orthogonal polynomials $\kappa_m\pi_m(z)$, $m\geq n$, (see \cite{BDJ}, (1.25)).


\section{Case $n=0$}\label{sec:0}
In this section, we prove \eqref{eq16} when $n=0$. 
From the Proposition \ref{lem1}, it is enough to show that 
\begin{equation}\label{eq78}
   \phi^{(1)}_{0}(\lambda)+ 
\biggl( -\frac{\partial}{\partial t}\biggr) \biggr|_{t=1} 
\det(I-\sqrt{t}\K_0) 
= \biggr[ \sum_{s=0}^{\infty}
\frac{\lambda^{s}}{(s!)^2} a_{0}(s)\biggr] \phi^{(1)}_{0}(\lambda).
\end{equation}
But this follows immediately from \eqref{eq17} 
in the case $n=0$, and the following Lemma.

\begin{lem}
  We have
\begin{equation}\label{eq79}
   1+ \frac12\tr\biggl( \frac1{I-\K_0}\K_0 \biggr)
= \sum_{s=0}^\infty \frac{\lambda^s}{(s!)^2} a_0(s).
\end{equation}
\end{lem}

\begin{proof}
From \eqref{e+18}, \eqref{e+19} and \eqref{eq40}, we have
\begin{equation}
  \begin{split}
     &tr\biggl( \frac1{I-\K_0}\K_0 \biggr)
= \int_\Sigma \bigl(F'(z;0;1)\bigr)^TG(z;0;1) dz \\
&\quad = \int_{|z|=1-\epsilon} 
\bigl(M(z;0;1)^{-1}M'(z;0;1)f(z;0;1)\bigr)^Tg(z;0;1) dz
+\int_{|z|=1-\epsilon} \bigl( f'(z;0;1)\bigr)^Tg(z;0;1)  dz,
  \end{split}
\end{equation}
for any $0<\epsilon <1$.
From the definition of $f, g$ with $k=0$ and $t=1$ in \eqref{e+10},
the second integral is $0$.
We use the relation between $M$ and $m$ in \eqref{e+7},
and the relation between $m$ and $Y$ in \eqref{e+6},  
to express the first integral in terms of $Y(z;0)$. Then we have
\begin{equation}\label{e+92}
   \tr\biggl( \frac1{I-\K_0}\K_0 \biggr)
= \int_{|z|=1-\epsilon}
\biggl( \frac12e^{\sqrt\lambda z}, -e^{-\sqrt\lambda z^{-1}}\biggr)
(Y(z;0)^{-1}Y'(z;0))^T
\biggl( 2e^{-\sqrt\lambda z}, e^{\sqrt\lambda z^{-1}} \biggr)^T \frac{dz}{2\pi i}.
\end{equation}
When $n=0$, a simple computation shows that 
\begin{equation}\label{eq-73}
   Y(z;0)= \begin{pmatrix}
1&\int_{\Sigma}\frac{\psi(s)}{s-z} \frac{ds}{2\pi i}\\ 0&1 \end{pmatrix}.
\end{equation}
Thus \eqref{e+92} becomes
\begin{equation}
  \begin{split}
   \tr\biggl( \frac1{I-\K_0}\K_0 \biggr)
&= -2\int_{|z|=1-\epsilon} \psi(z)^{-1}\frac{dz}{2\pi i}
\int_{\Sigma}\frac{\psi(s)}{(s-z)^2} \frac{ds}{2\pi i} \\
&= -2\int_{|z|=1-\epsilon} \frac{dz}{2\pi i}
\int_{\Sigma}\frac{e^{\sqrt\lambda(s-z)(1-\frac1{sz})}}{(s-z)^2} 
\frac{ds}{2\pi i}.
  \end{split}
\end{equation}
Now by Taylor expansion,
\begin{equation}
   \frac{e^{\sqrt\lambda(s-z)(1-\frac1{sz})}}{(s-z)^2}
= \sum_{k=0}^\infty \frac{(\sqrt\lambda)^k}{k!} 
(s-z)^{k-2}\bigl(1-\frac1{sz}\bigr)^k
= \sum_{k=0}^\infty \frac{(\sqrt\lambda)^k}{k!} \sum_{l=0}^{k}
\binom{k}{l}(s-z)^{k-2} \frac{(-1)^l}{s^lz^l}.
\end{equation}
For $|z|=1-\epsilon$, we have
\begin{equation}
  \int_{\Sigma}\frac{e^{\sqrt\lambda(s-z)(1-\frac1{sz})}}{(s-z)^2} 
\frac{ds}{2\pi i}
= \sqrt\lambda -\sum_{k=2}^\infty \frac{(-\sqrt\lambda)^k}{k!} 
\sum_{l=1}^{k-1}
\binom{k}{l}\binom{k-2}{l-1} z^{k-2l-1}.
\end{equation}
Therefore
\begin{equation}
  \tr\biggl( \frac1{I-\K_0}\K_0 \biggr)
= 2\sum_{k\geq 2, even} \frac{(-\sqrt\lambda)^k}{k!}
\binom{k}{k/2}\binom{k-2}{k/2-1}
= 2\sum_{p=1}^\infty \frac{\lambda^p}{(2p)!} \binom{2p}{p}\binom{2p-2}{p-1}.
\end{equation}
Thus we have
\begin{equation}\label{eq87}
    1+ \frac12\tr\biggl( \frac1{I-\K_0}\K_0 \biggr)
= 1+ \sum_{p=1}^\infty \frac{\lambda^p}{(p!)^2} \binom{2p-2}{p-1}.
\end{equation}

On the other hand, from the definition of $a_0(s)$ in \eqref{e+35}, 
the right hand side of \eqref{eq79} is given by 
\begin{equation}
  \sum_{s=0}^\infty \frac{\lambda^s}{(s!)^2} a_0(s)
= 1+ \sum_{s=0}^\infty \frac{\lambda^s}{(s!)^2}
\sum_{m=1}^{\infty} \frac{\lambda^{m}}{(m+s)^2((m-1)!)^2}.
\end{equation}
Changing the summation index to $p:= s+m$, this expression becomes 
\begin{equation}
  1+ \sum_{p=1}^\infty \sum_{s=0}^{p-1} \binom{p-1}{s}^2 \frac{\lambda^p}{(p!)^2}
= 1 + \sum_{p=1}^\infty \binom{2p-2}{p-1} \frac{\lambda^p}{(p!)^2}, 
\end{equation}
which agrees with \eqref{eq87}. 
Here we have used the elementary combinatorial identity 
\begin{equation}
   \sum_{s=0}^{p-1} \binom{p-1}{s}^2
= \binom{2p-2}{p-1}.
\end{equation}
This proves the Lemma.
\end{proof}




\section{Case $n>0$}\label{sec:positive}

From \eqref{eq49} and \eqref{eq52}, we need to verify that 
\begin{equation}
1-\sum_{k=n}^{\infty} \frac{\dot{m}_{11}(0;k+1;1)}{m_{11}(0;k+1;1)}
= \sum_{s=0}^{\infty}
\frac{\lambda^{s}}{(s!)^2} \biggl( a_{0}(s)-
\bigl( T^{-1}_{n-1}b_{n}(s),b_{n}(s)\bigr)\biggr).
\end{equation}
But using \eqref{eq49} and \eqref{eq78}, it is enough to show
\begin{equation}\label{eq80}
   -\sum_{k=0}^{n-1} \frac{\dot{m}_{11}(0;k+1;1)}{m_{11}(0;k+1;1)}
= \sum_{s=0}^{\infty}
\frac{\lambda^{s}}{(s!)^2} 
\bigl( T^{-1}_{n-1}b_{n}(s),b_{n}(s)\bigr).
\end{equation}
To this end, we need the following two Lemmas.
Recall $\psi(z) := e^{\sqrt\lambda(z+z^{-1})}$.

\begin{lem}\label{lem2}
   Let $0<\epsilon <1$.
   For $p,q\geq 0$, we have
\begin{equation}\label{eq93}
   \sum_{s=0}^{\infty} \frac{\lambda^s}{(s!)^2} a_{p+1}(s)a_{q+1}(s) 
= -\int_{|z|=1-\epsilon} \psi(-z)
  \biggl( \int_{|u|=1} \frac{u^q\psi(u)}{u-z} \frac{du}{2\pi i} \biggr)
  \biggl( \int_{|v|=1} \frac{\psi(v)}{v^{p+1}(v-z)} \frac{dv}{2\pi i} 
\biggr) \frac{dz}{2\pi i}.
\end{equation}
\end{lem}

\begin{proof}
  From the definition of $a_k(s)$, the left hand side of \eqref{eq93} is 
 \begin{equation*}
 \begin{split}
   &\sum_{s=0}^{\infty} \frac{\lambda^s}{(s!)^2} 
\sum_{a=p+1}^{\infty} \frac{\lambda^{a-(p+1)/2}}{(a+s)(a-1)!(a-p-1)!}
\sum_{b=q+1}^{\infty} \frac{\lambda^{b-(q+1)/2}}{(b+s)(b-1)!(b-q-1)!}\\
  &= \lambda^{-\frac{p+q}{2}} \sum_{s=0}^{\infty}\sum_{a=p+1}^{\infty}
\sum_{b=q+1}^{\infty} \frac{\lambda^{s+a+b-1}}{(s+a)!(s+b)!(a-p-1)!(b-q-1)!} 
\binom{s+a-1}{a-1} \binom{s+b-1}{b-1}.
 \end{split}
 \end{equation*}
 Setting $k=s+a+b-1-p-q$, $l=a-p$ and $m=b-q$, 
the above expression becomes
 \begin{equation}\label{e-LHS-1}
   \lambda^{\frac{p+q}{2}} \sum_{k=1}^{\infty} \lambda^k
\sum_{l=1}^{k} \sum_{m=1}^{k-l+1} \frac1{(k-m+p+1)!(k-l+q+1)!(l-1)!(m-1)!}
\binom{k-m+p}{l+p-1} \binom{k-l+q}{m+q-1}.
 \end{equation}

On the other hand, the right hand side of \eqref{eq93} is 
 \begin{equation}\label{e-RHS-1}
 -\int_{|z|=1-\epsilon} \frac{dz}{2\pi i}
\int_{|u|=1} \frac{u^qe^{\sqrt\lambda u}e^{\sqrt\lambda \frac{(z-u)}{uz}}}{u-z}
\frac{du}{2\pi i}
\int_{|v|=1} \frac{e^{\sqrt\lambda v^{-1}}e^{\sqrt\lambda (v-z)}}{v^{p+1}(v-z)}
\frac{dv}{2\pi i}.
 \end{equation}
For the third integral in \eqref{e-RHS-1}, Taylor expansions of 
$e^{\sqrt\lambda (v-z)}$ and $e^{\sqrt\lambda v^{-1}}$ give us 
 \begin{equation*}
   \sum_{a=1}^{\infty} \frac{(\sqrt\lambda)^a}{a!} \sum_{b=0}^{\infty} 
\frac{(\sqrt\lambda)^b}{b!}
\int_{|v|=1} \frac{(v-z)^{a-1}}{v^{b+p+1}} \frac{dv}{2\pi i} , 
 \end{equation*}
since the integration when $a=0$ vanishes.
Evaluating the integral, we obtain 
\begin{equation}\label{e-RHS-2}
   \sum_{a=p+1}^{\infty} \frac{(\sqrt\lambda)^a}{a!} \sum_{b=0}^{a-p-1} 
\frac{(\sqrt\lambda)^b}{b!} \binom{a-1}{b+p}(-z)^{a-b-p-1}.
\end{equation}
Similarly by expanding  
$e^{\sqrt\lambda \frac{(z-u)}{uz}}$ and $e^{\sqrt\lambda u}$ in Taylor 
series, the second integral in \eqref{e-RHS-1} becomes 	
\begin{equation}\label{e-RHS-3}
   z^qe^{\sqrt\lambda z} + \sum_{c=q+1}^{\infty} \frac{(\sqrt\lambda)^c}{c!} 
\sum_{d=0}^{c-q-1}
\frac{(\sqrt\lambda)^d}{d!} \binom{c-1}{d+q}(-z)^{-c+d+q}.
\end{equation}
Now using \eqref{e-RHS-2} and \eqref{e-RHS-3}, \eqref{e-RHS-1} becomes  
\begin{equation}\label{e-RHS-4}
   \sum_{a=p+1}^{\infty} \sum_{b=0}^{a-p-1} \sum_{c=q+1}^{\infty} 
\sum_{d=0}^{c-q-1} \frac{(-1)^{q-p}(-\sqrt\lambda)^{(a+b+c+d)}}{a!b!c!d!} 
\binom{a-1}{b+p} \binom{c-1}{d+q} 
\int_{|z|=1-\epsilon} z^{a-b-c+d-p+q-1} \frac{dz}{2\pi i}
\end{equation}
since the contribution from $z^qe^{\sqrt\lambda z}$ is zero.
The integral in \eqref{e-RHS-4} is $1$ 
when $d=-a+b+c+p-q$, which gives a new restriction on $c$,
$c\geq a-b-p+q$, implying 
\begin{equation*}
   \sum_{a=p+1}^{\infty} \sum_{b=0}^{a-p-1} \sum_{c=a-b-p+q}^{\infty} 
\frac{\lambda^{(b+c+(p-q)/2)}}{a!b!c!(-a+b+c+p-q)!}     
\binom{a-1}{b+p} \binom{c-1}{-a+b+c+p}. 
\end{equation*}
Changing the summation indices to $k,l,m$ by $b=l-1$, $c+b-q=k$ and 
$-a+b+c+p-q+1=m$, the above expression is \eqref{e-LHS-1}.
\end{proof}


\begin{lem}\label{lem3}
  Let 
\begin{equation}\label{e-26}
   \pi_n(z)=z^n+\cdots=\sum_{p=0}^n \eta^{n}_pz^p, \qquad \eta^{n}_n=1
\end{equation}
be the $n$-th monic orthogonal polynomial 
with respect to a measure $f(e^{i\theta})d\theta/2\pi$ on the unit circle 
which satisfies $f(e^{i\theta})=f(e^{-i\theta})$.
Let $\T_n=(c_{j-k})_{0\le j,k\le n}$ denote the $(n+1)\times (n+1)$ 
Toeplitz matrix with respect to the same measure. 
Then we have for $n\geq 1$,
\begin{eqnarray}
    \bigl(\T_n^{-1}\bigl)_{pq} &= \kappa_n^2\eta^{n}_p\eta^{n}_q, \qquad
&p=n\ \ \text{or}\ \ q=n, \\
    \bigl(\T_n^{-1}\bigl)_{pq}-\bigl(\T_{n-1}^{-1}\bigl)_{pq} 
&= \kappa_n^2\eta^{n}_p\eta^{n}_q, \qquad &0\le p,q\le n-1,
\end{eqnarray}
where $\kappa_n$ is the leading coefficient of the $n$-th 
normalized orthogonal polynomial, 
$\int_0^{2\pi} |\kappa_n\pi_n(e^{i\theta})|^2f(e^{i\theta})d\theta/2\pi=1$.
\end{lem}

\begin{proof}
   Set 
\begin{equation}
   \gamma_{pq}:= 
\begin{cases}
 \bigl(\T_n^{-1}\bigl)_{pq},  &p=n\ \ \text{or}\ \ q=n, \\
 \bigl(\T_n^{-1}\bigl)_{pq}-\bigl(\T_{n-1}^{-1}\bigl)_{pq}, 
&0\le p,q\le n-1.
\end{cases}
\end{equation}
Define 
\begin{equation}
  b(z,w) := \sum_{p,q=0}^{n} \gamma_{pq}z^pw^{q}.
\end{equation}
Using the definition of the Toeplitz coefficient 
$c_k=\int_0^{2\pi} e^{-ik\theta}f(e^{i\theta})d\theta/2\pi$,  
\begin{equation}\label{e-31}
 \begin{split}
  \int_0^{2\pi} b(e^{i\theta},w)e^{-ij\theta} 
f(e^{i\theta})\frac{d\theta}{2\pi} 
&= \sum_{p,q=0}^n \gamma_{pq}c_{j-p}w^{q}\\
&= \sum_{p,q=0}^n\bigl(\T_n^{-1}\bigl)_{pq}c_{j-p}w^{q}
- \sum_{p,q=0}^{n-1} \bigl(\T_{n-1}^{-1}\bigl)_{pq}c_{j-p}w^{q}.
 \end{split}
\end{equation}
But $c_{j-p}=\bigl(\T_n\bigl)_{jp}$ for $0\le j,p\le n$ and 
also $c_{j-p}=\bigl(\T_{n-1}\bigl)_{jp}$ for $0\le j,p\le n-1$.
Therefore \eqref{e-31} is zero for $0\le j\le n-1$.
This shows that for fixed $w$, $b(z,w)$ is a polynomial in $z$ 
of degree $n$ which is orthogonal to $1,z,\cdots,z^{n-1}$. 
Thus, for some $a(w)$, 
\begin{equation}
  b(z,w) = \pi_n(z)a(w).
\end{equation}
Now the evenness of $f$, $f(e^{i\theta})=f(e^{-i\theta})$, implies that 
$c_k=c_{-k}$, and hence 
the Toeplitz matrices above are symmetric, $\gamma_{pq}=\gamma_{qp}$, 
and so $b(z,w)=b(w,z)$. Thus
\begin{equation}
  b(z,w) = c\pi_n(z)\pi_n(w),
\end{equation}
for some constant $c$.
To determine the constant $c$, we consider the coefficient of the leading 
term $z^nw^n$ of $b(z,w)$.
That is 
\begin{equation}
  \gamma_{nn}=\bigl(\T_n^{-1}\bigl)_{nn}
=\frac{\det(\T_{n-1})}{\det (\T_n)}=\kappa_n^2 ,
\end{equation} 
by \cite{Sz}. 
Thus we have 
\begin{equation}
  b(z,w) = \kappa_n^2\pi_n(z)\pi_n(w)
\end{equation}
and this completes the proof.
\end{proof}

Finally we prove \eqref{eq80} which in turn completes 
the proof of Theorem \ref{thm-1}.
Let $m(z;k;t)$ be the solution of the RHP $(v(z;t),\Sigma)$ 
given in \eqref{e+5}. 

\begin{lem}\label{lem4}
We have for all $n\geq 1$,
\begin{equation}
   -\sum_{k=0}^{n-1}\frac{\dot{m}_{11}(0;k+1;1)}{m_{11}(0;k+1;1)}
=\sum_{s=0}^{\infty} \frac{\lambda^s}{(s!)^2} 
\bigl(\T_{n-1}^{-1}b_n(s),b_n(s)\bigr),
\end{equation}
\end{lem}

\begin{proof}
In the proof that follows, 
$m$,$\dot{m}$,$v$,$\dot{v}$ are all evaluated at $t=1$.
By differentiating the RHP \eqref{e+5} with respect to $t$, we have 
\begin{equation}
  \begin{cases}
    &\dot{m}_+=\dot{m}_-v+m_-\dot{v}, \qquad \text{on $\Sigma$},\\
    &\dot{m}=O(1/z) \qquad \text{as $z\to\infty$}.
  \end{cases}
\end{equation}
Since $m$ satisfies $m_+=m_-v$, we have 
\begin{equation}
  \begin{cases}
    &\bigl(\dot{m}m^{-1}\bigr)_+
=\bigl(\dot{m}m^{-1}\bigr)_- + m_-\dot{v}v^{-1}m_-^{-1}
=\bigl(\dot{m}m^{-1}\bigr)_- + m_+v^{-1}\dot{v}m_+^{-1}, 
\qquad \text{on $\Sigma$},\\
    &\dot{m}m^{-1}=O(1/z) \qquad \text{as $z\to\infty$}.
  \end{cases}
\end{equation}
From the Plemelj formula, the solution of this equation is given by 
\begin{equation}
   \bigl(\dot{m}m^{-1}\bigr)(z) = \int_\Sigma 
\frac{\bigl(m_+v^{-1}\dot{v}m_+^{-1}\bigr)(s)}{s-z} \frac{ds}{2\pi i}.
\end{equation}
Therefore for any $0<\epsilon <1$, we have 
\begin{equation}\label{e-42}
   \dot{m}_{11}(0) = \biggl[ \int_{|z|=1-\epsilon} 
m(z)v^{-1}(z)\dot{v}(z)m^{-1}(z)m(0) 
\frac{dz}{2\pi iz} \biggr]_{11}.
\end{equation}
To simplify the calculation note that the symmetry of the jump matrix 
$\bigl(\begin{smallmatrix} 0&1\\1&0\end{smallmatrix}\bigr) v(1/z) 
\bigl(\begin{smallmatrix} 0&1\\1&0\end{smallmatrix}\bigr) = v(z)^{-1}$
implies the symmetry of the solution 
\begin{equation}
    m(0)^{-1}m(1/z) = \begin{pmatrix} 0&1\\1&0\end{pmatrix}
m(z) \begin{pmatrix} 0&1\\1&0\end{pmatrix}.
\end{equation}
Using this relation for $m^{-1}(z)m(0)$ in \eqref{e-42}, we have 
\begin{equation}\label{e-44}
   \dot{m}_{11}(0) = \biggl[ \int_{|z|=1-\epsilon} 
m(z)v^{-1}(z)\dot{v}(z)\begin{pmatrix} 0&1\\1&0\end{pmatrix}m^{-1}(1/z)
\begin{pmatrix} 0&1\\1&0\end{pmatrix} \frac{dz}{2\pi iz} \biggr]_{11}.
\end{equation}
Note that $\dot{v}(z;k+1;1)=\biggl( 
\begin{smallmatrix} -1&-1/2z^{-k-1}\varphi^{-1}\\
1/2z^{k+1}\varphi&0 \end{smallmatrix} \biggr)$.
Now we express \eqref{e-44} in terms of $Y$ using the relation \eqref{e+6}, 
\begin{equation}\label{e-46}
 \begin{split}
   \dot{m}_{11}(0) = 
&\int_{|z|=1-\epsilon} \biggl[\frac12\psi(z)Y_{21}(z)Y_{21}(1/z) 
+\frac12z^{-k-1}Y_{21}(z)Y_{22}(1/z) \\
&\qquad\qquad   - \frac12z^{k+1}Y_{22}(z)Y_{21}(1/z) 
-\psi(z)^{-1}Y_{22}(z)Y_{22}(1/z) \biggr]\frac{dz}{2\pi iz},
 \end{split}
\end{equation}
where $\psi(z) = e^{\sqrt\lambda(z+z^{-1})}$ as given in \eqref{e+1}.

Now we use the explicit expression of $Y$ in terms of orthogonal polynomial 
given in \eqref{e+3}.
Especially we use the following expressions,  
\begin{eqnarray}
    Y_{21}(z) &=& -\kappa^2_{k}z^k\pi_{k}(1/z), \\
    Y_{21}(1/z) &=& -\kappa^2_{k}z^{-k}\pi_{k}(z), \\
    Y_{22}(z)  &=& -\kappa^2_{k} \int_{|v|=1}\frac{\pi_{k}(1/v)\psi(v)}{v(v-z)}
\frac{dv}{2\pi i},\\
    Y_{22}(1/z)  &=& \kappa^2_{k} \int_{|u|=1}\frac{z\pi_{k}(u)\psi(u)}{u-z}
\frac{du}{2\pi i}.
\end{eqnarray}
The first two terms in \eqref{e-46} cancel each other, which can be 
seen as follows : 
\begin{equation}
 \begin{split}
   \text{first term} &= \frac{\kappa^4_{k}}{2}\int_{|z|=1-\epsilon} 
\psi(z)\pi_{k}(1/z)\pi_{k}(z)\frac{dz}{2\pi iz},\\
   \text{second term} &= \frac{-\kappa^4_{k}}{2}\int_{|z|=1-\epsilon} 
\pi_{k}(1/z)
\biggl( \int_{|u|=1}\frac{\pi_{k}(u)\psi(u)}{u-z}
\frac{du}{2\pi i} \biggr) \frac{dz}{2\pi iz} \\
&= \frac{-\kappa^4_{k}}{2} \int_{|u|=1}\pi_{k}(u)\psi(u) \frac{du}{2\pi i}
\int_{|z|=1-\epsilon} \frac{\pi_{k}(1/z)}{u-z} \frac{dz}{2\pi iz} \\
&= \frac{-\kappa^4_{k}}{2} \int_{|u|=1}\pi_{k}(u)\psi(u) 
\pi_k(1/u)\frac{du}{2\pi iu} , 
 \end{split}
\end{equation}
which is $-1$ times the first term, by Cauchy.
Also the third term in \eqref{e-46} vanishes as 
\begin{equation}
  \begin{split}
    \text{third term} &= \frac{-\kappa^4_{k}}{2} \int_{|z|=1-\epsilon}
\pi_k(z) \biggl( \int_{|v|=1}\frac{\pi_{k}(1/v)\psi(v)}{v(v-z)}
\frac{dv}{2\pi i} \biggr) \frac{dz}{2\pi i} \\
&= \frac{-\kappa^4_{k}}{2} \int_{|v|=1} \pi_{k}(1/v)\psi(v) \frac{dv}{2\pi iv}
\int_{|z|=1-\epsilon} \frac{\pi_k(z)}{v-z}\frac{dz}{2\pi i} =0,
  \end{split}
\end{equation}
since the quantity $\frac{\pi_k(z)}{v-z}$ is analytic for $|z|\le 1-\epsilon$.
Thus together with the fact $m_{11}(0;k+1;1)=\kappa^2_{k}$, 
which follows immediately from \eqref{e+3} and \eqref{e+6}, we have
\begin{equation}\label{e-54}
   -\frac{\dot{m}_{11}(0;k+1;1)}{m_{11}(0;k+1;1)} 
= -\kappa^2_{k} \int_{|z|=1-\epsilon} 
\psi(z)^{-1} \biggl( 
\int_{|u|=1}\frac{\pi_{k}(u)\psi(u)}{u-z} \frac{du}{2\pi i} \biggr)
\biggl( \int_{|v|=1}\frac{\pi_{k}(1/v)\psi(v)}{v(v-z)} \frac{dv}{2\pi i} 
\biggr) \frac{dz}{2\pi i}.
\end{equation}

As in \eqref{e-26}, let 
\begin{equation}
   \pi_k(u)=\sum_{q=0}^{k} \eta^{k}_q u^q, \qquad 
\pi_k(1/v)=\sum_{p=0}^{k} \eta^{k}_p v^{-p}.
\end{equation}
Then \eqref{e-54} becomes 
\begin{equation}
   -\frac{\dot{m}_{11}(0;k+1;1)}{m_{11}(0;k+1;1)} 
= -\kappa^2_{k} \sum_{p,q=0}^{k}
\eta^{k}_q\eta^{k}_p \int_{|z|=1-\epsilon}
\psi(z)^{-1} \biggl( 
\int_{|u|=1}\frac{u^q\psi(u)}{u-z} \frac{du}{2\pi i} \biggr)
\biggl( \int_{|v|=1}\frac{\psi(v)}{v^{p+1}(v-z)} \frac{dv}{2\pi i} 
\biggr) \frac{dz}{2\pi i},
\end{equation}
and hence by Lemma \ref{lem2}, 
\begin{equation}\label{eq126}
   -\frac{\dot{m}_{11}(0;k+1;1)}{m_{11}(0;k+1;1)} = 
\sum_{s=0}^{\infty} \frac{\lambda^s}{(s!)^2} \sum_{p,q=0}^{k}
\kappa^2_{k} \eta^{k}_p\eta^{k}_q a_{p+1}(s)a_{q+1}(s).
\end{equation}
But then from Lemma \ref{lem3}, we have for $k\geq 1$,
\begin{equation}
   -\frac{\dot{m}_{11}(0;k+1;1)}{m_{11}(0;k+1;1)} = 
\sum_{s=0}^{\infty} \frac{\lambda^s}{(s!)^2} \biggl[ \sum_{p,q=0}^{k}
\bigl(\T_k^{-1}\bigr)_{pq} a_{p+1}(s)a_{q+1}(s) - \sum_{p,q=0}^{k-1}
\bigl(\T_{k-1}^{-1}\bigr)_{pq} a_{p+1}(s)a_{q+1}(s) \biggr].
\end{equation}
For $k=0$, $\eta^0_0= 1$, $\T_0$ is the $1\times 1$ matrix with 
entry $\kappa_0^{-2}$, and 
so by \eqref{eq126}, 
\begin{equation}
   -\frac{\dot{m}_{11}(0;1;1)}{m_{11}(0;1;1)} =
\sum_{s=0}^{\infty} \frac{\lambda^2}{(s!)^2} 
\bigl(\T_0^{-1}\bigr) a_{1}(s)a_{1}(s).
\end{equation}
Thus we have 
\begin{equation}
   -\sum_{k=0}^{n-1}\frac{\dot{m}_{11}(0;k+1;1)}{m_{11}(0;k+1;1)}
=\sum_{s=0}^{\infty} \frac{\lambda^s}{(s!)^2}
\sum_{p,q=0}^{n-1} \bigl(\T_{n-1}^{-1}\bigr)_{pq} a_{p+1}(s)a_{q+1}(s)
=\sum_{s=0}^{\infty} \frac{\lambda^s}{(s!)^2}
\bigl(\T_{n-1}^{-1}b_n(s),b_n(s)\bigr).
\end{equation}
\end{proof}


\section{Asymptotics}\label{sec:asymp}

In this section, we make some remarks concerning the proofs of 
Theorem \ref{thm-2} and \ref{thm-3}.

First, as in \cite{Jo}, one can show that 
$q^{(2)}_{n,N}$ is monotonically decreasing 
in $N$, and so by the de-Poissonization Lemma (see Lemma 2.5 in \cite{Jo}), 
it is enough to control $\phi^{(2)}_n(\lambda)$ as $n,\lambda\to\infty$. 
By \eqref{eq16}, this translates into controlling $\det(I-\sqrt{t}\K_n)$ 
for values of $t$ near $1$ (of course, the asymptotic behavior of 
$\phi^{(1)}_n(\lambda)$ is given in \cite{BDJ}). 
But as noted in the Introduction, $\K_n$, and hence $\sqrt{t}\K_n$, 
is an integrable operator, to which there is an canonically 
associated RHP (see \eqref{e+7.5}).
As in \cite{BDJ}, the steepest descent method can be used to 
analyze this RHP asymptotically as $\lambda,n\to\infty$.
Again the critical region is where $n\sim 2\sqrt\lambda$, 
in which case the RHP localizes to a small neighborhood of $z=-1$. 
Write $2\sqrt\lambda=k-xk^{1/3}/2^{1/3}$, where $x$ lies in a bounded set.
Writing $z=-1+s$ for $z$ near $-1$, we obtain 
\begin{equation}
 \begin{split}
  v(z;k;t) &= 
\begin{pmatrix}
1-t&-\sqrt{t}z^{-k}e^{-\sqrt\lambda(z-z^{-1})}\\
\sqrt{t}z^{k}e^{\sqrt\lambda(z-z^{-1})}&1 \end{pmatrix} \\
&= \begin{pmatrix}
1-t&-\sqrt{t}z^{-k}(-1)^ke^{-h(s,k,x)}\\
\sqrt{t}z^{k}(-1)^ke^{h(s,k,x)}&1 \end{pmatrix}, 
 \end{split}
\end{equation}
where
\begin{equation}
 \begin{split}
  h(s,k,x) &= -\frac{x}{2^{1/3}}(k^{1/3}s)
-\frac{x}{2^{4/3}k^{1/3}}(k^{1/3}s)^2 
+\frac{(k^{1/3}s)^3}6\bigl(1-\frac{3x}{k^{2/3}2^{1/3}}\bigr) +\cdots\\
&\sim -\frac{x}{2^{1/3}}(k^{1/3}s) +\frac{(k^{1/3}s)^3}6, 
\qquad \text{as $k\to\infty$}.
 \end{split}
\end{equation}
Rescaling $\tilde{s}=k^{1/3}s/2^{4/3}$, we see that we are lead to 
a RHP with jump matrix 
\begin{equation}
  \tilde{v} = 
\begin{pmatrix}
1-t&-(-1)^k\sqrt{t}e^{-2(-x\tilde{s}+\frac43\tilde{s}^3)}\\
(-1)^k\sqrt{t}e^{2(-x\tilde{s}+\frac43\tilde{s}^3)}&1 \end{pmatrix}
\end{equation}
on the line $i\R$ (cf. Figure 9 in \cite{BDJ}). 
But after rotating by $\pi/2$, this is precisely the RHP for the 
Painlev\'e II equation with parameters $p=-q=\sqrt{t}$, $r=0$ 
(cf. Figure 4 in \cite{BDJ} : the terms $(-1)^k$ can be removed by 
a simple conjugation).
These parameters $p,q,r$ correspond to the solution $u(x;t)$ 
of the Painlev\'e II equation, 
$u_{xx}=2u^3+xu$, with the boundary condition $u(x;t)\sim -\sqrt{t}Ai(x)$ 
as $x\to +\infty$, where $Ai$ is the Airy function (cf. \cite{BDJ} (1.4)).
As in Lemmas 5.1 and 6.3 in \cite{BDJ}, we can obtain an expression 
for $m_{11}(0;k+1;t)$ in terms of the solution of the above 
Painlev\'e II RHP. 
Inserting this information into \eqref{eq-49} in Proposition \ref{prop-1}, 
we learn that for $2\sqrt\lambda=n-xn^{1/3}/2^{1/3}$, 
$(1+\sqrt{t})^{-n}\det(I-\sqrt{t}\K_n) \to F(x;t)$ as $n\to\infty$.
Substituting this relation into \eqref{eq16}, we obtain the proof 
of Theorem \ref{thm-2}, and eventually Theorem \ref{thm-3}.
\medskip

\centerline{-  -  -  -  -  -  -  -  -  -  -  -  -  -  -  -  -
  -  -  -  -  -  -  -  -  -  -  -  -  -  -  -  -  -  -  -
  -  -  -  -  -  -  -  -  -  -  -  -  -  -  -  -}

\medskip
We conclude with some remarks on the motivation for a formula such as 
\eqref{eq16}. 
We started backwards, assuming that the second row behaves statistically 
in the large $N$ limit like the second largest eigenvalue of 
a GUE random matrix. 
As noted in the Introduction, this conjecture was strongly supported 
by numerical simulations of Odlyzko and Rains.
We had to end up with the Tracy-Widom distribution $F^{(2)}(x)$. 
From the point of view of \cite{BDJ}, $F^{(2)}(x)$ would have to 
emerge from the solution of some local RHP. 
For the case $F(x;1)$, which is expressed (see \eqref{e-Ft}) 
in terms of a specific solution $u(x;1)$ of the Painlev\'e II equation, 
$u_{xx}=2u^3+xu$, $u(x;1)\sim -Ai(x)$ as $x\to +\infty$, 
the local RHP was precisely the RHP for the Painlev\'e II equation, 
and this local problem emerged naturally via the steepest descent method 
applied to the RHP associated to Gessel's 
formula \eqref{eq13} in the canonical way.
But now $F^{(2)}$ is a derivative of $F(x;t)$ where $F(x;t)$ 
involves a family of solutions $\{ u(x;t) \}$ of the Painlev\'e II equation, 
$u_{xx}=2u^3+xu$, $u(x;t)\sim -\sqrt{t}Ai(x)$ as $x\to +\infty$. 
So we need to find a RHP which reduces in the critical region 
$n\sim 2\sqrt\lambda$ to a local RHP, which is precisely the RHP for 
the solution $u(x;t)$ of the Painlev\'e II equation. 
The RHP \eqref{e+5} is chosen precisely to ensure this property.

The procedure leading to \eqref{eq16} is now forced. 
The RHP \eqref{e+5} (more precisely, the equivalent RHP \eqref{e+7.5}) 
is of the type that arises from an integrable operator 
($\K_n$ in this case), which then leads after some calculations to the 
determinant formula for $I-\K_n$ on the right hand side of \eqref{eq16}.


\section*{Appendix}

In this Appendix, we first discuss the spectral properties of the operator 
$\K_n$ in \eqref{e-16} that are used in the proof of 
Lemma \ref{lem-1}, and then the
(unique) solvability of the RHP \eqref{e+5}. 

Let $\Sigma$ denote the unit circle in the complex plane, 
oriented counterclockwise, and let 
$\varphi(z)=e^{\sqrt\lambda(z-z^{-1})}$, as before.
The operator $\K_n : L^2(\Sigma, |dz|) \to L^2(\Sigma, |dz|)$ 
is defined by 
\begin{equation}\label{ap1}
   \K_n(z,w) := \frac{z^{-n}w^n-
\varphi(z)\varphi(w)^{-1}}{2\pi i(z-w)},
\qquad   \bigl(\K_n f\bigr)(z) = \int_\Sigma \K_n(z,w)f(w) dw.
\end{equation}
First note that 
\begin{equation}\label{ap2}
  \K_n= \frac12 A\widehat{H}B, 
\end{equation}
where the operators $A : L^2(\Sigma, |dz|)\oplus L^2(\Sigma, |dz|) 
\to L^2(\Sigma, |dz|)$ 
and $B : L^2(\Sigma, |dz|) \to 
L^2(\Sigma, |dz|)\oplus L^2(\Sigma, |dz|)$ 
are defined by 
\begin{equation}\label{ap3}
   \bigl(A\vec{h}\bigr)(z):=z^{-n}h_1(z)+\varphi(z)h_2(z), 
\qquad \bigl(Bh\bigr)(z):= 
\bigl( -z^nh(z), \varphi(z)^{-1}h(z) \bigr)^T,
\end{equation}
for a scalar $h$ and a vector $\vec{h}=(h_1, h_2)^T$,
and $\widehat{H} : L^2(\Sigma, |dz|)\oplus L^2(\Sigma, |dz|) 
\to L^2(\Sigma, |dz|)\oplus L^2(\Sigma, |dz|)$ is defined by 
\begin{equation}\label{ap4}
   \bigl(\widehat{H}\vec{h}\bigr)(z) 
:= \bigl( (Hh_1)(z), (Hh_2)(z)\bigr)^T, 
\end{equation}
for a vector $\vec{h}=(h_1, h_2)^T$ where 
$H : L^2(\Sigma, |dz|) \to L^2(\Sigma, |dz|)$ 
is the Hilbert transformation given by 
\begin{equation}\label{ap5}
   \bigl( Hh\bigr) (z) = \lim_{\epsilon\to 0} \frac1{i\pi} 
\int_{|s|=1, |s-z|>\epsilon} \frac{h(s)}{s-z} ds.
\end{equation}
Since $\|A\|\le \sqrt{2}$, $\|Bh\|=\sqrt{2}\|h\|$ and $\|Hh\|=\|h\|$, 
we have $\|\K_n \| \le 1$.

As $\varphi(z)=\overline{\varphi(z)^{-1}}$, we have 
\begin{equation}\label{ap6}
   \K_n(z,w)dw = 
\frac{\overline{z^n}w^n-\overline{\varphi(z)^{-1}}\varphi(w)^{-1}}
{2\pi\bigl(\overline{z^{-1}}w^{-1}-1 \bigr)} \frac{dw}{iw}. 
\end{equation}
Since $\frac{dw}{iw} = d\theta = |dw|$, $\K_N$ is a self-adjoint operator 
on $L^2(\Sigma, |dz|)$.
Also since the kernel is smooth, the operator $\K_n$ is trace class, 
and hence $\|\K_n\|=1$ if and only if $+1$ and/or $-1$ is 
an eigenvalue.

We show that $1$ is not an eigenvalue of $\K_n$. 
Observe first that 
\begin{equation}\label{ap7}
   \|A\vec{h}\|= \sqrt{2}\|\vec{h}\| \quad 
\text{if and only if} \quad z^{-n}h_1(z)=\varphi(z)h_2(z).
\end{equation}
Now suppose that $1$ is an eigenvalue of $\K_n$. 
Then there is a non-trivial function $h\in L^2(\Sigma, |dz|)$ 
such that $\K_n h = h$.
Then 
\begin{equation}\label{ap8}
   \|h\| = \|K_n h\| = \frac12 \|A\widehat{H}Bh\| 
\le \frac1{\sqrt{2}} \|\widehat{H}Bh\| = \|h\|,
\end{equation}
which implies that 
\begin{equation}\label{ap9}
   \|A\widehat{H}Bh\| = \sqrt{2} \|\widehat{H}Bh\|.
\end{equation}
Hence by \eqref{ap7} above, and by the definition of 
the operator $B$ given in \eqref{ap3}, we have an equation
\begin{equation}\label{ap10}
   -z^{-n}H(z^nh) = \varphi(z)H(\varphi(z)^{-1}h).
\end{equation}
Now re-express $\K_n$ as follows,
\begin{equation}\label{ap11}
   \K_nh= -\frac12z^{-n}H(z^nh) + 
\frac12\varphi(z)H(\varphi(z)^{-1}h).
\end{equation}
Using \eqref{ap10}, 
\begin{equation}\label{ap12}
   \K_nh= -z^{-n}H(z^nh) = \varphi(z)H(\varphi(z)^{-1}h),
\end{equation}
which leads to the equations 
\begin{equation}\label{ap13}
   -z^{-n}H(z^nh) = h, \qquad 
\varphi(z)H(\varphi(z)^{-1}h) = h,
\end{equation}
or 
\begin{equation}\label{ap14}
   H(z^nh) = -z^{n}h, \qquad 
H(\varphi(z)^{-1}h) = \varphi(z)^{-1}h.
\end{equation}
But $H(z^n)=z^n$ for $n\geq 0$ 
and $H(z^n)=-z^{n}$ for $n<0$,
so that \eqref{ap14} implies 
\begin{equation}\label{ap15}
   z^nh(z)=\sum_{j<0}a_jz^j, \qquad 
\varphi(z)^{-1}h(z) = \sum_{j\geq 0}b_j z^j, 
\end{equation}
for some square summable sequences 
$\{a_j\}_{j<0}$ and $\{b_j\}_{j\geq 0}$.
Hence, the second equation in \eqref{ap15} implies that 
\begin{equation}\label{ap16}
   e^{\sqrt\lambda z^{-1}}h(z) = e^{\sqrt\lambda z} \sum_{j\geq 0}b_j z^j.
\end{equation}
Combining with the first equation in \eqref{ap15}, 
\begin{equation}\label{ap17}
   e^{\sqrt\lambda z^{-1}}\sum_{j<0}a_jz^{j-n}
= e^{\sqrt\lambda z} \sum_{j\geq 0}b_j z^j,
\end{equation}
which is impossible for $n\geq 0$ 
unless all the $a_j$'s (and $b_j$'s) are zero.
Therefore for $n\geq 0$, $1$ is not an eigenvalue of $\K_n$.
In a similar manner, one can show that 
$\dim Ker(\K_n+1) =n$, for $n\geq 0$.
In particular, $\|\K_n\|=1$ and $\dim Ker(\K_n-1)=0$ 
for $n\geq 0$.

Now we prove the (unique) solvability of the RHP \eqref{e+5}.
It is clear that the solvability of the RHP \eqref{e+5} 
follows from the solvability of the RHP \eqref{e+7.5} since 
$m$ and $M$ are algebraically related by \eqref{e+7}. 
Now from integrable operator theory (see Lemma 2.21 \cite{DIZ}), 
the existence of the inverse of $(I-\sqrt{t}\K_k)^{-1}$ 
implies the solvability of the RHP \eqref{e+7.5}.
But as $\|K_k\|=1$, and as $1$ is not in the spectrum of $\K_k$, 
it follows that $(I-\sqrt{t}\K_k)^{-1}$ exit for all $0<t\le 1$, 
and hence the RHP \eqref{e+5} is solvable.
The proof of the uniqueness of the solution of the RHP is standard 
(compare, for example, \cite{BDJ} Lemma 4.2). 

\noindent {\bf Remark.} Note from \eqref{eq17} that when $\lambda=0$, 
$\det(I-\K_n(\lambda=0))=2^n\phi^{(1)}_n(0)=2^n$,
which can be checked directly (for $\lambda=0$, $-\K_n(\lambda=0)$ 
is an orthogonal projection of rank $n$).
On the other hand, as $\lambda, n\to\infty$, the spectrum of 
$\K_n=\K_n(\lambda)$ can approach $1$ : indeed from Lemma 7.1 (v) 
in \cite{BDJ}, for $2\sqrt\lambda\geq (n+1)(1+\delta_7) \to\infty$, 
$\delta_7>0$, $\phi^{(1)}_n(\lambda)\le Ce^{-cn^2}$ so that 
$\det(I-\K_n)=2^n\phi^{(1)}_n(\lambda)\to 0$ as $n\to\infty$.


\end{document}